\documentclass[12pt]{article}
\usepackage{hyperref}
\usepackage{amsfonts,amsmath,amstext,amssymb}
\usepackage{dsfont}
\usepackage{graphicx}
\usepackage{multicol}
\usepackage{mathtools}
\usepackage{mathptmx}

\newtheorem{thm}{Theorem}

\title{On the hypersufaces of the Euclidean space which are simultaneously minimal and maximal}

\author{
Magdalena Caballero\footnote{Partially supported by Spanish MINECO and FEDER Project MTM2016-78807-C2-1-P.}\\
{\tt \small magdalena.caballero@uco.es}\\
\and
{\small Departamento de Matem\'aticas, Campus de Rabanales,}\\ {\small Universidad de C\'ordoba, 14071 C\'ordoba, Spain}}

\date{}

\begin{document}

\maketitle

\begin{abstract}
	It is well known that the only surfaces that are simultaneously minimal in $\mathbb{R}^3$ and maximal in $\mathbb{L}^3$ are open pieces of helicoids (in the region in which they are spacelike)  and of spacelike planes, \cite{Kobayashi}. The proof of this result consists in showing that the level curves of those surfaces are lines, and so the surfaces are ruled. And it finishes comparing the classification of minimal ruled surfaces to that of maximal ruled surfaces.
	
	In this manuscript we consider the general case of spacelike hypersurfaces in the $(n+1)$-dimensional Euclidean space which are simultaneously maximal and minimal. We show that its level curves are minimal hypersurfaces in the $n$-dimensional Euclidean space.
\vspace{.3cm}

\noindent {\bf Keywords: spacelike hypersurfaces, maximal hypersurfaces, minimal hypersurfaces} 

\noindent {\bf 2010 MSC: 53C42; 53C50; 35J93} 
\end{abstract}

\section{Introduction and background}

A hypersurface in the Lorentz-Minkowski space $\mathbb{L}^{n+1}$ is said to be spacelike if its induced metric is a Riemannian one. Spacelike hypersurfaces in $\mathbb{L}^{n+1}$ can be endowed with another Riemannian metric, the one inherited from the Euclidean space $\mathbb{R}^{n+1}$. Consequently, we can consider two different mean curvature functions on a spacelike hypersurface, $H_R$ and $H_L$ respectively. 

It is well known that any spacelike hypersurface $\Sigma$ in $\mathbb{L}^{n+1}$ can be locally described as a spacelike graph over an open subset of a spacelike hyperplane, which without loss of generality can be supposed to be the hyperplane $x_{n+1}=0$ (see~\cite{Lopez} for the two dimensional case). 

Let $u$ be a function over an open subset of $\mathbb{R}^n$ and let $\Sigma_u$ be its graph. Then the spatiality condition is written as $|Du|<1$. It is possible to get expressions for the normal vector field to $\Sigma_u$ in $\mathbb{R}^{n+1}$, $N_R$, and in $\mathbb{L}^{n+1}$, $N_L$, as well as for the mean curvature functions $H_L$ and $H_R$, in terms of $u$. Specifically, with a straightforward computation we get
$$\label{eq:normLR}N_R=\dfrac{(-Du,1)}{\sqrt{1+|Du|^2}} \,\textrm{,} \qquad N_L=\dfrac{(Du,1)}{\sqrt{1- |Du|^2}},$$ 
$$
\cos\theta=\frac{1}{\sqrt{1+|Du|^2}} \,\textrm{,} \qquad \cosh\psi=\frac{1}{\sqrt{1-|Du|^2}},
$$
$$
H_R = \dfrac{1}{n} \, \mathrm{div}\left(\dfrac{Du}{\sqrt{1+|Du|^2}}\right) \,\textrm{,}\qquad \textrm{and} \qquad
H_L = \dfrac{1}{n} \, \mathrm{div}\left(\dfrac{Du}{\sqrt{1- |Du|^2}}\right), 
$$
where $D$, $\mathrm{div}$ and $|\cdot|$ stand for the gradient, the divergence and the Euclidean norm on $\mathbb{R}^n$, respectively, whereas $\theta$ and $\psi$ denote the angle between $N_R$ and $(0,...,0,1)$ and the hyperbolic angle between $N_L$ and $(0,...,0,1)$, respectively.

Therefore, a spacelike graph determined by $u$ satisfies $H_R=H_L$ if and only if $u$ is a solution of the following partial differential equation 
$$
\mathrm{div}\left(\left(\frac{1}{\sqrt{1-|Du|^2}}- \frac{1}{\sqrt{1+|Du|^2}}\right)Du\right)=0,
\text{ with } |Du|<1,
$$
called the {\em $H_R=H_L$ hypersurface equation}. This equation is a quasilinear elliptic partial differential equation, everywhere except at those points at which $Du$ vanishes, where the equation is parabolic, see~\cite{AlarconAlbujerCaballero}.

As a particular case, we can consider the situation where the graph is simultaneously minimal and maximal, that is $H_R=H_L=0$. The Calabi-Bernstein theorem, proved by Calabi for dimension $3$, see \cite{Calabi}, and by Cheng and Yau for arbitrary dimension, see \cite{ChengYau}, states that the only entire maximal graphs in $\mathbb{L}^{n+1}$ are the spacelike hyperplanes. As an immediate consequence we conclude that the only entire graphs that are simultaneously minimal in $\mathbb{R}^{n+1}$ and maximal in $\mathbb{L}^{n+1}$ are the spacelike hyperplanes. 

Going a step further, we can consider spacelike graphs with the same constant mean curvature functions. Heinz~\cite{Heinz}, Chern~\cite{Chern} and Flanders~\cite{Flanders} proved that the only entire graphs with constant mean curvature $H_R$ in $\mathbb{R}^{n+1}$ are the minimal graphs. Taking into account that any complete spacelike hypersurface in $\mathbb{L}^{n+1}$ is necessarily an entire graph over any spacelike hyperplane, see~\cite[Proposition 3.3]{AliasRomeroSanchez}, and using again the Calabi-Bernstein theorem, we conclude that the only complete spacelike hypersurfaces in $\mathbb{L}^{n+1}$ with the same constant mean curvature functions $H_R$ and $H_L$ are the spacelike hyperplanes.

Kobayashi~\cite{Kobayashi} studied the surfaces with $H_R=H_L=0$ without assuming any global hypothesis. He showed that the graphs which are simultaneously minimal and maximal are open pieces of a spacelike plane or of a helicoid, in the region where the helicoid is a spacelike surface. In general dimension, Lee and Lee~\cite{LeeLee} have recently presented non-planar examples of simultaneously minimal and maximal spacelike graphs in the Lorentz-Minkowski space. Their examples can be seen as generalized ruled surfaces, in fact they are a natural gene\-ralization of helicoids. However, there is no known classification of such hypersurfaces.

Recently, Albujer, Caballero and S\'anchez~\cite{AlbujerCaballero,AlbujerCaballeroSanchez} have continued with the study of spacelike surfaces with the same mean curvature in $\mathbb{R}^3$ and in $\mathbb{L}^3$, not necessarily constant. The generalization of some of its results to general dimension has been given in \cite{AlarconAlbujerCaballero, AlarconAlbujerCaballero2}. 

The goal of this work is to prove the following new result. 
\begin{thm}
	Let us consider a hypersurface in $\mathbb{R}^{n+1}$ which is simultaneously minimal and maximal. Then its level hypersurfaces are minimal surfaces in $\mathbb{R}^n$.
\end{thm}

This theorem is a corollary of \cite[Lema 4.1]{AlarconAlbujerCaballero2} which was not included in \cite{AlarconAlbujerCaballero2}. In this manuscript we are going to include all the computations without taking into account the previously cited lemma, giving a slightly different proof.

\section{Proof of the theorem}

As we have mentioned before, any spacelike hypersurface is locally a graph over the $x_{n+1}\equiv 0$ hyperplane, which can be identified with $\mathbb{R}^n$. Let $u$ be a function over an open set $\Omega\subseteq \mathbb{R}^n$ whose graph, $\Sigma_u$, is simultaneously minimal and maximal and let $\pi:\Sigma_u\longrightarrow \Omega$ denote the canonical projection.

Given $p\in\Sigma_u$ such that $Du(\pi(p))\neq 0$, we consider its corresponding level hypersurface contained in $\mathbb{R}^n$, $\widetilde{S_c}$, and its lifting to $\Sigma_u$, $S_c$. We will work in a neighborhood of $p$ on which $Du\circ\pi$ does not vanish. Hence the  distribution given by $Du$ is integrable, so we can consider the integral curve through $\pi(p)$. We denote by  $\alpha$ its lifting to $\Sigma_u$. Notice that $ \alpha' = (Du, |Du|^2) \circ \pi $. 

Therefore, we have two submanifolds of $\Sigma_u$, namely $S_c$ and $\alpha$, defined on a neighborhood of $p$ which are orthogonal at $p$ for both $\langle\cdot,\cdot\rangle_R$ and $\langle\cdot,\cdot\rangle_L$, the metrics in $\mathbb{R}^{n+1}$ and $\mathbb{L}^{n+1}$, respectively. Now, let  $ \{e_1, \ldots, e_{n-1}\} $ be an orthonormal basis of $ T_{\pi(p)} \widetilde{S_c}$. The vectors $ \{(e_1,0), \ldots, (e_{n-1},0)\} $ constitute an orthonormal basis of $ T_p S_c$ in both $\mathbb{R}^{n+1}$ and $\mathbb{L}^{n+1}$, and are orthogonal to $\alpha'$ for both metrics. 

Then, using the definition of the normal curvature along a direction of a hypersurface in $\mathbb{R}^{n+1}$ and in $\mathbb{L}^{n+1}$, respectively, as well as \eqref{eq:normLR}, we obtain the following relationships, where we have omitted the point $ p $ on behalf of simplicity, and where $\kappa_v^R$ is the normal curvature of $\Sigma_u$ at $p$ along the direction $v$ in $\mathbb{R}^n$, while in $\mathbb{L}^n$ it is denoted by $\kappa_v^L$.

\[
\kappa_{(e_i,0)}^R= - \dfrac{\cos\theta}{\cosh\psi}\,\kappa_{(e_i,0)}^L 
, \qquad i = 1, \ldots, n - 1, \qquad \text{and}
\] 

\[
\kappa_{\alpha'}^R= -\dfrac{\cos^3\theta}{\cosh^3\psi}\,\kappa_{\alpha'}^L. 
\] 

By denoting $A=\dfrac{\cos\theta}{\cosh\psi}$, we rewrite the previous relationships as 
\begin{equation}\label{eq:relkA1}
\kappa_{(e_i,0)}^R=-A\,\kappa_{(e_i,0)}^L, \,\, i = 1, \ldots, n - 1, \,\, \text{and} \,\,\, \kappa_{\alpha'}^R=-A^3\,\kappa_{\alpha'}^L.
\end{equation}

As we are dealing with orthogonal directions at $p$ for both $\langle\cdot,\cdot\rangle_R$ and $\langle\cdot,\cdot\rangle_L$, and $u$ is a solution of the    
$H_R=H_L$ hypersurface equation, we get
\begin{equation*}\label{eq:relkA2}
-\kappa_{(e_1,0)}^L-\ldots-\kappa_{(e_{n-1},0)}^L-\kappa_{\alpha'}^L=\kappa_{(e_1,0)}^R+\ldots+\kappa_{(e_{n-1},0)}^R+\kappa_{\alpha'}^R,
\end{equation*}
which jointly with \eqref{eq:relkA1} implies 
$$-\kappa_{\alpha'}^L=\dfrac{1}{A^2+A+1}\left(\kappa_{(e_1,0)}^L+\ldots+\kappa_{(e_{n-1},0)}^L\right).
$$

But since the hypersurface is maximal, we also know that 
$$-\kappa_{\alpha'}^L=\kappa_{(e_1,0)}^L+\ldots+\kappa_{(e_{n-1},0)}^L.
$$

Since $Du(\pi(p))\neq 0$, $A\neq 1$, and so we conclude that
\begin{equation}\label{relacion-kLalfa-KLbeta2}-\kappa_{\alpha'}^L=\kappa_{(e_1,0)}^L+\ldots+\kappa_{(e_{n-1},0)}^L=0. 
\end{equation}

On the other hand, for each $i=1,\dots,n-1$ we take a curve in $\widetilde{S_c}$, $\widetilde{\alpha_i}$, with $\widetilde{\alpha_i}(0)=p$ and $\widetilde{\alpha_i}'(0)=e_i$. Let $\alpha_i$ be its lifting to $S_c$. Notice that $\alpha_i'=(\widetilde{\alpha}_i',0)$. 

Let us observe that it is possible to relate the Lorentzian normal curvature $\kappa_{(e_i,0)}^L$ of $\Sigma_u$ at $p$ in the direction of $(e_i,0)$ with the normal curvature $\kappa_{e_i}^c$ of $\widetilde{S_c}$ at $\pi(p)$ in the direction of $e_i$:

\[
\kappa_{(e_i,0)}^L=\langle\overline{\nabla}_{t_i}t_i,N_L\rangle_L=\frac{|Du|\left\langle D_{\widetilde{t_i}}\,\widetilde{t_i},\dfrac{Du}{|Du|}\right\rangle_{\mathbb{R}^n}}{\sqrt{1-|Du|^2}}\circ\pi=\frac{|Du|\kappa_{e_i}^c}{\sqrt{1-|Du|^2}}\circ\pi.
\]
Here $\overline{\nabla}$ is the Levi-Civita connection in $\mathbb{L}^{n+1}$, $D$ and $\langle\cdot,\cdot \rangle_{\mathbb{R}^n}$ stand for the Levi-Civita connection and the usual metric of the Euclidean space $\mathbb{R}^n$, respectively, $t_i=\dfrac{\alpha'_i}{|\alpha'_i|_L}$,  $\widetilde{t_i}=\dfrac{\widetilde{\alpha_i}'}{|\widetilde{\alpha_i}'|}$ and $\dfrac{Du}{|Du|}$ is the unitary normal vector field to $\widetilde{S_c}$ in $\mathbb{R}^n$. 

Therefore, from~\eqref{relacion-kLalfa-KLbeta2} we get that the mean curvature of $\widetilde{S_c}$ in $\mathbb{R}^n$ vanishes at $p$. 

By a continuity argument, the result can be proved true for points at which the gradient vanishes.

\bibliographystyle{amsplain}

\end{document}